\documentclass[11pt]{article}          
\usepackage{amsmath,amssymb} 

\newtheorem{theorem}{Theorem}           
\newtheorem{lemma}{Lemma}

\newtheorem{remark}{Remark}[section]

\newcommand{\proof}{\noindent {\it Proof}. }
\newcommand{\eop}{\hfill $\Box$}

\newcommand{\xx}{\mathbf x}
\newcommand{\yy}[1]{\mathbf{y}^{(#1)}}
\newcommand{\Zsum}[2]{S_{#1,#2}}
\newcommand{\upbnd}{C^+}  
\begin{document}

\title%
{Lower bound for cyclic sums of Diananda type}

\author{Sergey Sadov\footnote{%
E-mail: serge.sadov@gmail.com}}

\date{}                               

\maketitle

\begin{abstract}
Let $C=\inf (k/n)\sum_{i=1}^n x_i(x_{i+1}+\dots+x_{i+k})^{-1}$, where the infimum is taken over all pairs of integers
$n\geq k\geq 1$ and all positive $x_1,\dots,x_{n+k}$
subject to cyclicity assumption $x_{n+i}=x_i$, $i=1,\dots,k$. We prove that $\ln 2\leq C< 0.9305$. In the definition of the constant $C$ the operation $\inf_k\inf_n\inf_{\xx}$ can be replaced by $\lim_{k\to\infty}\lim_{n\to\infty}\inf_{\xx}$.


\end{abstract}


\section{Background and main theorem}
Given integers $n\geq k\geq 1$ and a vector $\xx$ with positive components 
 $x_1,\dots,x_n$ (we write $\xx>0$),
let us define  {\it interval sums}
$$
 t_{i,k}=\sum_{j=0}^{k-1} x_{i+j}
$$
and the {\it cyclic sum of Diananda type}
$$
 \Zsum{n}{k}(\xx)=\sum_{i=1}^n \frac{x_i}{t_{i+1,k}}.
$$
 Hereinafter we treat subscripts modulo $n$, that is,   
 $x_{n+i}=x_i$ by definition. 
  
For example,
$$
 \Zsum{3}{2}(\xx)=\frac{x_1}{x_2+x_3}+\frac{x_2}{x_3+x_1}+\frac{x_3}{x_1+x_2}.
$$

The sums $\Zsum{n}{2}$ are commonly referred to as Shapiro sums. Following \cite[p.~217]{BoDay_1973}, we associate the name of 
P.H.~Diananda with the more general sums $\Zsum{n}{k}$ 
because they have first occured in his note \cite{Diananda_1959}.   
   
The function $\Zsum{n}{k}(\xx)$ is homogeneous of degree zero in its vector argument.
Let 
\begin{equation}
\label{Ank}
 A(n,k)=\inf_{\xx>0}\, \Zsum{n}{k}(\xx).
\end{equation}
The domain of the function  $\Zsum{n}{k}(\xx)$ can be extended to allow zero values of some $x_i$:
it suffices to require that $t_{i,k}>0$ for all $i$. The value of $\inf$ is not affected. 
 
For every $k=1,2,\dots$ denote
\begin{equation}
\label{Bk}
 B(k)=\inf_{n\geq k}\, \frac{k}{n} A(n,k).
\end{equation}%
Define
\begin{equation}
\label{Ck}
C=\inf_{k\geq 1} B(k).
\end{equation}

\begin{theorem}
\label{thm:main}
Let $\upbnd\approx 0.930498$ be the $y$-intercept of the common tangent to the graphs
$y=e^{-x}$ and $y=x/(e^x-1)$. 
Then
\begin{equation}
\label{mainest}
\ln 2\leq C\leq \upbnd.
\end{equation}
\end{theorem}

This is our main result. Let us put it in context.

If all $x_i$ are equal, then
$$
 \Zsum{n}{k}(\xx)=\frac{n}{k},
$$
hence always
$$
A(n,k)\leq \frac{n}{k}, \qquad
B(k)\leq 1.
$$

For $k=1$ we have
$$
 \Zsum{n}{1}(\xx)=\frac{x_1}{x_2}+\frac{x_2}{x_3}+\dots+\frac{x_{n-1}}{x_n}+\frac{x_n}{x_1}\geq n
$$
by the inequality between arithmetic and geometric means (AM-GM). 
Thus $B(1)=1$.

The inequality $\Zsum{3}{2}(\xx)\geq 3/2$ implying
$A(3,2)= 3/2$  has been known since 1903 at latest
(due to A.M.~Nesbitt) --- see Ref.\ \cite[p.~440]{MPF_1993}, which offers three proofs.
H.S.~Shapiro \cite{Shapiro_1954} proposed to prove that $A(n,2)\geq n/2$ for all $n$, i.e.\ that $B(2)=1$. This conjecture 
was soon disproved \cite{Northover_1956}.  The precise validity range for Shapiro's conjectured inequality
was determined through analytical and numerical labor over time span of more than 20 years: even $n\leq 12$ and odd $n\leq 23$.
See review \cite{Clausing_1992}.  
The actual value of $B(2)$  found by  V.G.~Drinfeld \cite{Drinfeld_1971} is slightly less than one: $B(2)=0.989133\dots$. It equals the $y$-intercept of the common tangent to the graphs $y=g_1(x)=e^{-x}$ and $y=g_2(x)=2/(e^{x/2}+e^x)$. 
 
Much less was known until now about lower bounds for the cyclic sums $\Zsum{n}{k}$
with $k\geq 3$. 
In \cite{Diananda_1961} Diananda showed that
$$
\frac{k}{n}A(n,k)\geq \frac{2(k+1)}{n}
$$ 
for $n>2(k+1)$. (In the 
region $2(k+1)<n<(2/\ln 2)(k+1)$ of the $(n,k)$-lattice this lower bound beats the estimate $(k/n)A(n,k)\geq \ln 2$ that follows from our Theorem~\ref{thm:main}.)
\@\ 
Diananda also found a few cases where $A(n,k)=n/k$ with $k>2$. 
They are listed with references in \cite[p.~173]{Diananda_1977} or \cite[p.~445]{MPF_1993}. 
%
%
The cited results do not allow one to conclude that $B(k)>0$ for $k\geq 3$. The only result to that effect was Diananda's \cite{Diananda_1962} simple lower bound $A(n,k)\geq n/k^2$, which implies that $B(k)\geq 1/k$. 
Compare: Theorem~\ref{thm:main} says that in fact $B(k)$ are uniformly separated from zero.
  
A systematic study of cyclic sums in which numerators and denominators are {\it overlapping}\ interval sums 
is carried out in \cite{Baston_1973}. The closest in appearance to our inequality
$A(n,k)\geq \mathrm{const}\cdot(n/k)$, $\mathrm{const}=\ln 2$, is 
Baston's formula (in our notation)
$$
\inf_{\xx} \,\sum_{i=1}\frac{x_i}{t_{i,k}}= \frac{1}{1+\left\lfloor\frac{k-1}{n}\right\rfloor}
\qquad (k\geq 2)
$$
contained in his Theorem~1. 
Yet the presence of the summand $1$ in the denominator on the right causes a striking contrast with our situation:
the analog of $B(k)$, 
$$
\inf_{n:\,n\geq k}\;\,\inf_{\xx} \,\frac{k}{n}\sum_{i=1}\frac{x_i}{t_{i,k}}, 
$$ 
equals zero unless $k=1$!

Some further relevant citations can be found in Remarks to Theorems~2--4 below.
For a detailed review of similar and other cyclic inequalities see \cite[Ch.~16]{MPF_1993}, particularly \S~15 and further on.

In the following two sections we will establish 
$k$-dependent bounds for the individual constants $B(k)$ tighter than their common bounds in Theorem~\ref{thm:main}; the latter will easlily follow. 
The lower and upper estimates are treated separately, since the methods of proof are different.

The final section of the paper answers in the affirmative a natural question whether the operations $\inf_{n:\,n\geq k}$ 
and $\inf_{k}$ in the definitions \eqref{Bk}
of $B(k)$ and  \eqref{Ck} of $C$ can be replaced, respectively, by $\lim_{n\to\infty}$ and $\lim_{k\to\infty}$.  

\section{Lower bounds for $B(k)$}
\label{sec:lbnd}

\begin{theorem}
\label{thm:lbnd}
The constants $B(k)$ are bounded below as follows:
$$
 B(k) \geq  k(2^{1/k}-1).
$$
In other words, for any integers $n\geq k\geq 1$ and any $n$-dimensional vector $\xx>0$
the cyclic inequality
$$
\frac{k}{n}\,\sum_{i=1}^n \frac{x_i}{t_{i+1,k}}\geq k(2^{1/k}-1)
$$
holds.
\end{theorem}

\begin{remark}\rm
The left inequality \eqref{mainest} of Theorem~\ref{thm:main} follows since $k(2^{1/k}-1)>\ln 2$ (indeed, $e^x-1>x$ ($x\neq 0$); take $x=k^{-1}\ln 2$).
\end{remark}

\begin{remark}\rm
The numerical values of our lower bounds for $k\leq 7$ are listed below.

\bigskip 
\noindent
\begin{tabular}{c|cccccc} 
$k$            & $2$ & $3$ & $4$ & $5$ & $6$ & $7$ \\ 
\hline
$k(2^{1/k}-1)$ & $0.82843$ & $0.77976$ & $0.75683$ & $0.74349$ & $0.73477$ & $0.72863$ 
\end{tabular}

\bigskip
In the case $k=2$ our value is worse than the best lower estimate $B(2)\geq 0.922476\dots$ \cite{Diananda_1962ii} known before Drinfeld's exact result
$B(2)=0.989133\dots$, yet it is better than earlier attempts, e.g. 
$B(2)\geq 0.66046\dots$ due to R.A.~Rankin \cite{Rankin_1961}.
This will be helpful to keep in mind when trying to improve our result.

The author is unaware of any published lower bound for $B(3)$ except for $B(3)\geq 1/3$, a particular case of Diananda's inequality
$B(k)\geq 1/k$.
\end{remark}

\proof
Without loss of generality we may assume that $k|n$ ($k$ divides $n$). Indeed, given any $k$ and $n$, 
let $n'=kn$ and define the $n'$-dimensional vector $\xx'$ as concatenation of $k$ copies of $\xx$. 
Obviously, $k|n'$ and
$$
\frac{1}{n'}\sum_{i=1}^{n'}\frac{x'_i}{t'_{i+1,k}}=\frac{1}{n}\sum_{i=1}^{n}\frac{x_i}{t_{i+1,k}}.
$$
(The index arithmetic for $\xx'$ is done modulo $n'$, unlike for $\xx$.)  
 
\smallskip 
From now on we assume that $n=k\nu$ with $\nu\geq 1$ an integer.
We fix the vector argument $\xx$ for the rest of the proof.
Set
$$
 r_j=\frac{t_{k j+1,k}}{t_{k (j+1)+1,k}}, \quad j=0,\dots, \nu-1.
$$
(When $j=\nu-1$ the denominator equals $t_{1,k}$.) 

We will derive lower estimates for the partial sums 
$$
 s_j=\sum_{i=1}^{k} \frac{x_{jk+i}}{t_{jk+i+1,k}}, \quad j=0,\dots, \nu-1.
$$
Consider $j=0$ for simplicity of notation.
We have
$$
x_i=t_{i,2k-i+1}-t_{i+1,2k-i}.
$$
Also for $i=1,\dots,k$ 
$$
 t_{i+1,k}\leq   t_{i+1,2k-i}.
$$
Hence
$$
 s_0\geq \sum_{i=1}^k \frac{t_{i,2k-i+1}-t_{i+1,2k-i}}{t_{i+1,2k-i}}
 =-k+\sum_{i=1}^k \frac{t_{i,2k-i+1}}{t_{i+1,2k-i}}.
$$
Since
$$
 \prod_{i=1}^k \frac{t_{i,{2k-i+1}}}{t_{i+1,2k-i}}=\frac{t_{1,2k}}{t_{k+1,k}}=\frac{t_{1,k}+t_{k+1,k}}{t_{k+1,k}}=1+r_0,
$$
the AM-GM inequality yields
$$
 s_0\geq k\left((1+r_0)^{1/k}-1\right).
$$
Similar inequalities hold for $s_j$ with $j=1,\dots,\nu-1$ in place of $j=0$.

Introduce the function
$$
  f_k(t)=k \left( (1+e^t)^{1/k}-1 \right).
$$
Let $r_j=e^{t_j}$. From the above we obtain
$$
\Zsum{n}{k}(\xx)=\sum_{j=0}^{\nu-1} s_j \geq \sum_{j=0}^{\nu-1} f_k(t_j). 
$$
The function  $f_k(t)$ is convex. A simple way to see it is to factor the derivative as follows
$$
 f_k'(t)=\frac{(1+e^t)^{1/k}}{1+e^{-t}},
$$
where the numerator increases while the denominator decreases,
so $f_k''(t)>0$.

Since $\prod_{j=0}^{\nu-1} r_j=1$, the Jensen inequality yields
$$
\frac{\Zsum{n}{k}(\xx)}{\nu}\geq f_k\left(\frac{1}{\nu}\sum_0^{\nu-1} t_j\right)=f_k(0)=k(2^{1/k}-1). 
$$
This is precisely the claimed inequality.
\eop

\section{Upper bounds for $B(k)$}
\label{sec:ubnd}

In this section we use Drinfeld's \cite{Drinfeld_1971} construction to prove upper estimates for $B(k)$, $k=2,3,\dots$. For $k=2$
we just get Drinfeld's constant.
It is conceivable 
that for every $k\geq 3$ this construction provides a minimizing sequence $\xx^{(n)}$, too, 
but we are unable to prove it. 
Drinfeld's proof (for $k=2$) does not extend to $k>2$.

As a preparation to Theorem~\ref{thm:ubnd} below let us introduce
the family of functions
\begin{equation}
\label{gk}
g_k(x)=\frac{k(1-e^{-x/k})}{e^x-1}, \qquad k\geq 1,
\end{equation}
and study some of their properties.
We set $g_k(0)=1$ to 
make $g_k(x)$ continuous (in fact, real analytic).
Note that $g_1(x)=e^{-x}$.
Denote also
$$
 g_\infty(x)=\lim_{k\to\infty} g_k(x)=\frac{x}{e^x-1}.
$$

\begin{lemma} 
\label{lem:ubnd1}
The functions $g_\infty(x)$, $p(x)=(1-e^{-x})/x$, and
$g_k(x)$ for any $k>0$ are positive, decreasing and convex.
\end{lemma}

\proof Note that $g_k(x)=g_\infty(x/k) p(x)$.
By the Leibniz Rule the class of positive, decreasing and convex functions is closed under multiplication. It is also  
closed under rescaling of the independent variable.
Hence it suffices to give a proof for $g_\infty(x)$ and $p(x)$.
The only nontrivial task is to check convexity.
For $g_\infty$ we have
$$
 g_\infty''(x) = 
 \frac{e^{2x}(x-2)+e^x(x+2)}{(e^x-1)^3}
 =\frac{\cosh \frac x2\,\left(\frac x2-\tanh \frac x2\right)}{2\sinh^3 \frac x2} >0.
$$ 
Now for $p(x)$: if $x<0$, then we write $p(x)=\tilde p(-x)$, where
$\tilde p(y)=(e^y-1)/y$ has the Maclaurin series with positive coefficients, hence convex when $y>0$. And if $x>0$,
then we calculate
$$
p''(x)=\frac{2e^{-x}}{x^3}\left(e^x-\left(1+x+\frac{x^2}{2}\right)\right),
$$
and conclude that $p''(x)>0$ by the Maclaurin expansion, again.
\eop

\begin{lemma}
\label{lem:ubnd2}
For a fixed real $x\neq 0$ the function $k\mapsto k(1-e^{-x/k})$ is increasing. 
\end{lemma}

\proof
The claim is immediately clear, for both signs of $x$, by writing
$$
k(1-e^{-x/k})=\int_0^x e^{-y/k}\,dy.
$$
\eop

\begin{lemma} 
\label{lem:ubnd3}
The function $k\mapsto g_k(x)$ is increasing for any fixed $x>0$
and decreasing for any fixed $x<0$. In particular, 
for $k>1$ the inequalities
$$
g_k(x)>g_1(x)=e^{-x}\quad \mbox{\rm if $x>0$}
$$
and
$$
g_k(x)<e^{-x}\quad \mbox{\rm if $x<0$}
$$
hold.
\end{lemma}

\proof
Apply Lemma~\ref{lem:ubnd1}, taking into account sign of
$e^x-1$.
\eop

\bigskip
Let $h_k(x)$ be the {\it convex minorant}\ of the function
$\min(g_1(x),g_k(x))$.
From Lemmas~\ref{lem:ubnd1} and \ref{lem:ubnd3} it follows that
$h_k(x)$ is of the form
$$
h_k(x)=\left\{
\begin{array}{ll}
 g_k(x), & x\leq a_k, \\[1ex]
 \gamma_k+\lambda_k x,  \; & a_k< x< b_k,
  \\[1ex]
 e^{-x}, & x\geq b_k,
\end{array}
\right.
$$
where $a_k<0<b_k$ are the abscissas of the tangency points
of the common tangent to the graphs $y=g_k(x)$ and $y=e^{-x}$.
The parameters $a_k$, $b_k$, $\gamma_k$ and $\lambda_k$ are uniquely determined by the condition
that $h_k(x)$ be continuous and differentiable.
A simple way to find them numerically is as follows (we omit the subscript $k$ to lighten notation). The tangent to the graph $y=e^{-x}$ at $(b,e^{-b})$ is $y=e^{-b}(1+b-x)$. It is also tangent to $y=g(x)$ at $(a,g(a))$, hence $-e^{-b}=g'(a)$ and
$g(a)=-g'(a)(1+b-a)$. Eliminating $b$ from the last two equations
leads to the equation for the single unknown $a$:
\begin{equation}
\label{comtan}
  \frac{g(a)}{g'(a)}-a+1=\ln(-g'(a)).
\end{equation}
Now $\lambda=g'(-a)$, $b=-\ln \lambda$, and $\gamma=-\lambda(1+b)$. 

Clearly, $\gamma_k<1$ and $\lambda_k<0$.
The pointwise monotonicity of the family 
$g_k(x)$ stated in Lemma~\ref{lem:ubnd3} implies that
$
\mbox{\rm as $k$ increases, $\gamma_k$ and $|\lambda_k|$ decrease}.
$
As $k\to\infty$, the parameters $a_k$, $b_k$, $\gamma_k$ and $\lambda_k$ tend to their limits corresponding to the 
convex minorant of $\min(g_1(x),g_\infty(x))$. Therefore the constant $\upbnd$ in Theorem~1 is the monotone limit
$$
 \upbnd =\gamma_\infty=\lim_{k\to\infty} \downarrow \gamma_k.
$$
Consequently, the upper bound in Theorem~\ref{thm:main}
(left inequality in \eqref{mainest}) follows from 
more precise estimates in the next theorem.
 
\begin{theorem}
\label{thm:ubnd}
Let $k\geq 2$, the function $g_k(x)$ be defined by \eqref{gk},
and $\gamma_k$ be (as defined above) the $y$-intercept of the common tangent to the graphs
$y=e^{-x}$ and $y=g_k(x)$.
The constant $B(k)$ defined in \eqref{Bk} satisfies the inequality 
$$
B(k)\leq \gamma_k.
$$ 
\end{theorem}   

\begin{remark}\rm
Some numerical values of our upper bounds (found by solving Eq.~\eqref{comtan}) are listed below.  The limit value $\gamma_\infty=\upbnd$ is included for convenience of comparison.

\bigskip 
\noindent
\begin{tabular}{c|cccccc|c}
$k$            & $2$ & $3$ & $4$ 
& $10$ 
& $100$ & $1000$ & $\infty$
\\
\hline
$\gamma_k$ & $0.98913$ & $0.97793$ & $0.96994$ 
& $0.94983$ &
$0.93272$ & $0.93072$ & $0.930498$
\end{tabular}

\bigskip
Here $\gamma_2$ is nothing but Drinfeld's constant.
Besides it, the only other previously reported estimate of this sort seems to
be that due to J.C.~Boarder and D.E.~Daykin \cite[Table 2, row `$a/bcd$']{BoDay_1973} (reproduced in \cite{MPF_1993} as Eq.~(27.41), p.~453):
$\inf_n A(n,3)/n\leq 0.32598$, implying $B(3)\leq 0.97794$.
Our estimate, with more digits than in the table above, is $B(3)\leq\gamma_3=0.9779277986\dots$. Since
$\gamma_3/3>0.32598-0.5\times 10^{-5}$, within the accuracy of 5 significant digits we can not claim an improvement over \cite{BoDay_1973}. 
However the method of  \cite{BoDay_1973} is entirely numerical and based on calculation of bounds for $A(n,3)/n$ for finitely many 
$n$, while in the proof below we let $n\to\infty$. 
\end{remark}

\proof
Fix $k\geq 2$. Given an $\epsilon>0$ we will find an integer $n>k$ and an $n$-dimensional vector $\xx$ such that
$\Zsum{n}{k}(\xx)<\gamma_k+\epsilon$.

The point $(0,\gamma_k)$ is a convex combination of the points $(a_k,g_k(a_k))$ and $(b_k,e^{-b_k})$ with some coefficients
$\mu_k$ and $\mu'_k=1-\mu_k$,
$$
 \mu_k a_k+\mu'_k b_k=0, \qquad \mu_k g_k(a_k)+\mu'_k e^{-b_k}=\gamma_k.
$$
Let us choose rational $\mu_*=p/q\in (0,1)$ sufficiently close to $\mu_k$
and real $a_*$, $b_*$ sufficiently close to $a_k$, $b_k$ respectively so that 
\begin{equation}
\label{mustar_eps}
 \mu_* a_*+\mu'_* b_*=0, \qquad \mu_* g_k(a_*)+\mu'_* e^{-b_*}<\gamma_k+\frac{\epsilon}{2}.
\end{equation}
From now on $\mu_*$, $a_*$, $b_*$ are assumed fixed.
We write $\mu_*$ as a fraction generally not in the lowest terms, 
$$
 \mu_*=\frac{m}{n}.
$$ 
Later it will be important to allow $n$ be as large as we please.
We may and will assume that the numerator $m$ and denominator $n$ are divisible by $k$.
   
Let us now describe construction of $\xx$ assuming $n$ and $m=\mu_* n$ given; a specific choice of $n$ will be made afterwards. 

Denote $m'=n-m=\mu'_* n$.   
Define the $n$-dimensional vector $\xx$ as follows:
$$
\begin{array}{ll}
 x_{jk}=e^{jb_*},\quad & \mbox{\rm if $1\leq j<m'/k$};
\\[1ex] %
 x_i=0,\quad &\mbox{\rm if $1\leq i<m'$, $k\not|\,i$};  
\\[1ex] %
 x_i=
 e^{-a_*(n-i)/k},\quad & \mbox{\rm if $m'\leq i\leq n$}.
\end{array}
$$ 

The sequence $x_i$ of length $n$ consists of two parts. It is sparse when $i<m'$; only one in every $k$ consecutive terms is nonzero, and those nonzero terms are increasing. 
For $m\leq i\leq n$, all terms are nonzero; they form a decreasing geomeric sequence.

Note that the formula $x_{jk}=e^{jb_*}$, which is the definition when $j<m'/k$, continues to hold for $j=m'/k$. 
%
%
Indeed, it follows from the equality $(m'/k) b_*  = -a_* (n-m')/k$, which is true since $\mu'_* b_*=- \mu_* a_*$.

Let us compute nonzero terms in the sum $\Zsum{n}{k}(\xx)$.
For $i=jk<m'$ we have
$$
\frac{x_i}{t_{i+1,k}}=\frac{x_i}{x_{i+k}}=e^{-b_*}.
$$
For $m'\leq i\leq n-k-1$ (there are $m-k$ such $i$'s) we get
$$
\frac{x_i}{t_{i+1,k}}=\left(\sum_{j=1}^k e^{a_* j/k}\right)^{-1}=
\frac{1-e^{-a_*/k}}{e^{a_*}-1}= \frac{g_k(a_*)}{k}.
$$
For the remaining $k$ terms (with $i=n-k,\dots,n-1$) a convenient closed form expression is not available.
The rough estimate $t_{i+1,k}> x_{i+1}=e^{a_*/k} x_i$ will suffice. Thus
$$
 \sum_{i=n-k}^{n-1}\frac{x_i}{t_{i+1,k}} < k e^{-a_*/k}.
$$

In total,
$$
 \Zsum{n}{k}(\xx)< \frac{m'}{k} e^{-b_*}+\frac{m-k}{k}g_k(a_*)+ ke^{-a_*/k}.
$$
So
$$
 \frac{k}{n} \Zsum{n}{k}(\xx)< \mu'_* e^{-b_*}+\mu_*  g_k(a_*)+\frac{\delta}{n},
$$
where
$$
 \delta=k^2 e^{-a_*/k}-k g_k(a_*)
$$
does not depend on $n$.  
We choose $n$ so as to make 
$$
 \frac{\delta}{n}<\frac{\epsilon}{2}.
$$
%
Recalling \eqref{mustar_eps}, we obtain
$$
\frac{k}{n}S_{n,k}(\xx)<\gamma_k+\epsilon.
$$
Since $\epsilon$ is  arbitrary, we conclude that
$B(k)\leq \gamma_k$.
\eop

\section{Greatest lower bounds $B(k)$ and $C$ as limits}

\begin{theorem}
\label{thm:lim}
The constants $B(k)$ and $C$ defined in \eqref{Bk} and \eqref{Ck}
respectively can be expressed as limits. Specifically,
\\[\medskipamount]
{\rm (a) }\ for every integer $k\geq 1$
$$
B(k)=\lim_{n\to\infty} \frac{k}{n} A(n,k);
$$
{\rm (b) }\ the sequence $B(k)$ is nonincreasing and (therefore)
$$
 C=\lim_{k\to\infty} \downarrow B(k).
$$
\end{theorem}

\begin{remark}\rm
It has been known since the early period of investigation of 
the original Shapiro's conjecture, that $A(n,2)/n$ is not monotone in $n$.
\end{remark}

\begin{remark}\rm
Part (a) of this Theorem was established for Shapiro sums ($k=2$)
by R.A.~Rankin \cite{Rankin_1958} and for very general cyclic sums --- by K.~Goldberg \cite{Goldberg_1960}. To make our exposition self-contained, we still give a proof, which is a generalization
(albeit obvious) of Rankin's and more explicit than Goldberg's.
\end{remark} 

\proof
(a)
We assume $k\geq 1$ fixed once for all.

Given an $\epsilon>0$, we will find $N$ such that 
$$
 \frac{A(n,k)}{n}< \frac{B(k)}{k}+\epsilon.
$$
whenever $n\geq N$.

There exist an $m>k$ and an $m$-dimensional vector $\xx$ for which 
$$
 \frac{\Zsum{m}{k}(\xx)}{m}< \frac{B(k)}{k}+\frac{\epsilon}{2}.
$$
Fix an arbitrary $(m-1)$-tuple $(\xi_1,\dots,\xi_{m-1})$
of positive numbers. 
Define the vectors (``$r$-extensions of $\xx$ by $\mathbf{\xi}$''):
$\yy{0}=\xx$, and for $r=1,\dots,m-1$ 
$$
\yy{r}=(x_1,\dots,x_m,\xi_1,\dots,\xi_r).
$$
Let
$$
 M=\max_{0\leq r\leq m-1}\Zsum{m+r}{k}(\yy{r})
$$
and choose $N$ big enough to make $M/N<\epsilon/2$.

Now, given $n=m\nu+r$, $\nu\geq 1$, $0\leq r\leq m-1$, we construct
an $n$-dimensional vector $\xx'$ as concatenation of $(\nu-1)$ copies of $\xx$ followed by $\yy{r}$. 
It is readily seen that
$$
\Zsum{n}{k}(\xx')=(\nu-1)\Zsum{m}{k}(\xx)+\Zsum{m+r}{k}(\yy{r}).
$$
Therefore, if $n\geq N$, then
$$
\frac{A(n,k)}{n}\leq \frac{\Zsum{n}{k}(\xx')}{n}<\frac{(\nu-1)m}{n}
\left(\frac{B(k)}{k}+\frac{\epsilon}{2}\right)+\frac{M}{n}
<\frac{B(k)}{k}+\frac{\epsilon}{2}+\frac{\epsilon}{2},
$$
as required.

\medskip
(b) 
Given a $k\nu$ dimensional vector $\xx$, define a $(k+1)\nu$ dimensional vector $\xx'$ as follows:
for $j=0,1,\dots,\nu-1$
$$
 x'_{(k+1)j+r}=\left\{\begin{array}{ll}
  x_{kj+r} \;& \mbox{\rm if $1\leq r\leq k$},
  \\[1ex]
  0 \;& \mbox{\rm if $r= k+1$}.
  \end{array}
 \right. 
$$
Then 
$$
 \Zsum{(k+1)\nu}{k+1}(\xx')=\Zsum{k\nu}{k}(\xx).
$$
Taking $\inf_{\xx}$ we conclude:
$$
 A((k+1)\nu,k+1)\leq A(k\nu,k).
$$
By part (a), 
$$
 \lim_{\nu\to\infty} \frac{A(k\nu,k)}{\nu }=B(k), 
$$
and the inequality $B(k+1)\leq B(k)$ follows.
\eop

\end{document}